# Base Dependence of Benford Random Variables

Frank Benford

February, 2017; revised April, 2021

## 1. Introduction.

My grandfather, the physicist Frank Benford for whom Benford's Law is named, considered his "law of anomalous numbers" as evidence of a "real world" phenomenon. He realized that geometric sequences and exponential functions are generally base 10 Benford, and on this basis he wrote (Benford, (1938) [1]):

> If the view is accepted that phenomena fall into geometric series, then it follows that the observed logarithmic relationship is not a result of the particular numerical system, with its base, 10, that we have elected to use. Any other base, such as 8, or 12, or 20, to select some of the numbers that have been suggested at various times, would lead to similar relationships; for the logarithmic scales of the new numerical system would be covered by equally spaced steps by the march of natural events. As has been pointed out before, the theory of anomalous numbers is really the theory of phenomena and events, and the *numbers* but play the poor part of lifeless symbols for living things.

This argument seems compelling, and it might seem to apply to Benford random variables as well as to geometric sequences and exponential functions. It is therefore somewhat surprising to observe that a random variable that is base $b$ Benford is not generally base $c$ Benford when $c \neq b$. We'll see some examples shortly.

This paper builds on two of my earlier papers ([2],[3]) and is an attempt to cast some light on the issue of base dependence. It's organized as follows. Section 2 introduces the significand function and the fractional part notation and gives several logically equivalent definitions of "Benford random variable." The base $b$ first digit law is introduced, and several examples of random variables are presented that are Benford relative to one base but not to another. Section 3 introduces the "Benford spectrum" $B_X$ of a positive random variable $X$ and summarizes some of the known analytical results that involve $B_X$. Section 4 is a brief digression listing some facts about Fourier transforms that are needed in subsequent sections. Section 5 introduces some fundamental notation and results that provide a framework for the "Benford analysis" of a positive random variable. Section 6 combines the framework of Section 5 with my method of "seed functions" to develop the theory of the base $c$ Benford properties of random variables $X$ that are known to be Benford relative to base $b$, and Section 7 gives several examples of such random variables. Section 8 discusses Berger and Hill's concept of "base-invariant significant digits." Section 9 is a summary and a look ahead.

## 2. Benford Random Variables

The best way to define Benford random variables is via the *significand function*. Let $b > 1$ be a fixed "base." Any $x > 0$ may be written uniquely in the form

$$x = s \times b^k$$

where $1 \leq s < b$ and $k \in \mathbb{Z}$, and the **base $b$ significand of $x$**, written $S_b(x)$, is defined as this $s$. Hence

$$x = S_b(x) \times b^k \tag{2.1}$$

where $1 \leq S_b(x) < b$ and $k \in \mathbb{Z}$. (Berger and Hill [4] define the significand of $x$ for all $x \in \mathbb{R}$, but we don't require this generality.)

Now let $X$ be a positive random variable; that is, $\Pr(X > 0) = 1$. Assume that $X$ is continuous with a probability density function (pdf). **Definition 2.1**: $X$ is base $b$ Benford (or $X$ is $b$-Benford) if and only if the cumulative distribution function (cdf) of $S_b(X)$ is given by

$$\Pr(S_b(X) \leq s) = \log_b(s) \quad \text{for all} \quad s \in [1, b). \tag{2.2}$$

Nothing written above requires that $b$ be an integer. For this paragraph *alone*, we assume that $b$ is an integer greater than or equal to 3. Let $D_1(X)$ denote the first (i.e., leftmost or most significant) digit of $X$ in the base $b$ representation of $X$, so $D_1(X) \in \{1, \ldots, b-1\}$. (Leading zeros, if there are any, are ignored.) **Proposition 2.1**. If $X$ is $b$-Benford, then

$$\Pr(D_1(X) = d) = \log_b\left(\frac{d+1}{d}\right) \tag{2.3}$$

for all $d \in \{1, \ldots, b-1\}$. This is the **base $b$ First Digit Law**. To prove it, it is sufficient to observe that $D_1(X) = d$ iff $d \leq S_b(X) < d + 1$.

It's useful at this point to introduce some non-standard notation. Let $y \in \mathbb{R}$ and recall that the "floor" of $y$, written $\lfloor y \rfloor$, is defined as the largest *integer* that is less than or equal to $y$. Define $\langle y \rangle$ as:

$$\langle y \rangle \equiv y - \lfloor y \rfloor \tag{2.4}$$

and note that $0 \leq \langle y \rangle < 1$ for every $y \in \mathbb{R}$. We'll call $\langle y \rangle$ the *fractional part* of $y$, though if $y < 0$ this description is misleading.

If we take the logarithm base $b$ of eq. (2.1) we obtain

$$\log_b(x) = \log_b(S_b(x)) + k. \tag{2.5}$$

On the other hand,

$$\log_b(x) = \langle \log_b(x) \rangle + \lfloor \log_b(x) \rfloor. \tag{2.6}$$

As $\lfloor \log_b(x) \rfloor$ is necessarily an integer and $0 \leq \log_b(S_b(x)) < 1$, comparison of eqs. (2.5) and (2.6) shows that

$$\log_b(S_b(x)) = \langle \log_b(x) \rangle \qquad \text{and} \qquad k = \lfloor \log_b(x) \rfloor \tag{2.7}$$

for any $x > 0$.

Using eq. (2.7), we may rephrase Definition 2.1 in several logically equivalent ways: $X$ is $b$-Benford iff

(1) $\Pr(\log_b(S_b(X)) \leq \log_b(s)) = \log_b(s)$ for all $s \in [1, b)$,
(2) $\Pr(\langle \log_b(X) \rangle \leq u) = u$ for every $u \in [0, 1)$,
(3) $\langle \log_b(X) \rangle \sim U[0, 1)$,
(4) $X = b^Y$ where $\langle Y \rangle \sim U[0, 1)$,

where the notation "$W \sim U[0, 1)$" means that $W$ is uniformly distributed on the half open interval $[0, 1)$. (More generally, I use the symbol " $\sim$ " to mean "is distributed as." Hence, for example, "$X \sim f()$" means that $X$ is distributed with pdf $f()$, "$X_1 \sim X_2$" means that $X_1$ and $X_2$ have the same distribution, and "$X \sim N(\mu, \sigma)$" means that $X$ is distributed normally with mean $\mu$ and variance $\sigma^2$.)

For any random variable $Y$, if $\langle Y \rangle \sim U[0, 1)$ we sometimes say that $Y$ is "uniformly distributed modulo one," abbreviated "u.d. mod 1." Hence $X$ is $b$-Benford iff $\log_b(X)$ is u.d. mod 1.

With this background we can now give a couple of examples of random variables that are Benford with respect to one base but not to another. Let $Y \sim U[0, 1)$. **Example 1**. Let $X \equiv 10^Y$, so $X$ is 10-Benford. But it's not 8-Benford as it fails to satisfy the base 8 First Digit Law. To see this, note that the support of $X$ is $[1, 10)$, and let $D_1(X)$ denote the first digit in the base 8 representation of $X$. Then

$$\begin{aligned} \Pr(D_1(X) = 1) &= \Pr(1 \leq X < 2) + \Pr(8 \leq X < 10) \\ &= \log_{10}(2) + \log_{10}(5/4) \approx 0.3979, \end{aligned}$$

whereas

$$\log_8\left(\frac{1+1}{1}\right) = \frac{\ln(2)}{\ln(8)} = \frac{1}{3}.$$

**Example 2**. Let $Y$ be as above, but now let $X \equiv 8^Y$, so $X$ is 8-Benford. Note that the support of $X$ is $[1, 8)$. Let $D_1(X)$ denote the first digit in the base 10 representation of $X$. Hence $\Pr(D_1(X) \in \{8, 9\}) = 0$, whereas

$$\log_{10}\left(\frac{9}{8}\right) + \log_{10}\left(\frac{10}{9}\right) \approx 0.09691.$$

Hence, $X$ fails to satisfy the base 10 First Digit Law.

## 3. The Benford Spectrum.

Let $X$ be a positive random variable. **Definition 3.1**. Following Wójcik [7], the "**Benford spectrum**" of $X$, denoted $B_X$, is defined as

$$B_X \equiv \{b \in (1, \infty)\colon X \text{ is } b\text{-Benford}\}. \tag{3.1}$$

The Benford spectrum of $X$ may be empty. In fact, the Benford spectra of essentially all the standard random variables used in statistics are empty.

This section summarizes some of the known facts about Benford spectra. While proofs are provided for Proposition 3.2 and 3.4, I'm just going to provide citations for proofs of the other propositions.

**Proposition 3.1** (Berger and Hill [4], page 44, Proposition 4.3 (iii)). A random variable $Y$ is u.d. mod 1 if and only if $kY + c$ is u.d. mod 1 for every integer $k \neq 0$ and every $c \in \mathbb{R}$.

The condition in this proposition that $k$ be an integer cannot be relaxed.

**Proposition 3.2** (Whittaker [6]). If $b \in B_X$, then $\sqrt[m]{b} \in B_X$ for all $m \in \mathbb{N}$. In other words, if $X$ is $b$-Benford, then $X$ is $\sqrt[m]{b}$-Benford for all $m \in \mathbb{N}$. **Proof**. Suppose that $X$ is $b$-Benford, so $X = b^Y$ where $Y$ is u.d. mod 1. Hence, for any $m \in \mathbb{N}$,

$$X = \left(b^{1/m}\right)^{mY}.$$

As $b^{1/m} > 1$ and $mY$ is u.d. mod 1 by Proposition 3.1, it follows that $b^{1/m} \in B_X$.

**Proposition 3.3**. If $B_X$ is non-empty, then it is bounded above. In other words, no random variable can be $b$-Benford for arbitrarily large $b$. **Citations**: [3], [6], [7].

**Proposition 3.4**. If $X$ is $b$-Benford and $c > 0$, then $cX$ is $b$-Benford. **Proof**. As $X$ is $b$-Benford, $Y \equiv \log_b(X)$ is u.d. mod 1. As $\log_b(cX) = Y + \log_b(c)$ is u.d. mod 1 from Proposition 3.1, it follows that $cX$ is $b$-Benford.

We say of this result that the Benford property of a random variable is "scale-invariant."

**Proposition 3.5**. Suppose that $X$ and $W$ are independent positive random variables and that $X$ is $b$-Benford. Then the product $XW$ is also $b$-Benford. **Citations**: [3], [4], [7].

**Proposition 3.6** (a corollary of Proposition 3.5). If $X$ and $W$ are independent positive random variables, then

$$B_X \cup B_W \subseteq B_{XY}.$$

So far, the spectra we've seen are at most countably infinite. One may wonder if there exists a random variable with an uncountable spectrum. Whittaker showed by an example that such a random variable exists. Let $b > 1$ be given. Define $g: \mathbb{R} \to \mathbb{R}$ by

$$g(y) \equiv \frac{1 - \cos(2\pi y)}{2\pi^2 y^2}. \tag{3.2}$$

It may be shown that $g()$ is a legitimate pdf, and $Y$ is u.d. mod 1 if $Y \sim g()$. Hence $X \equiv b^Y$ is $b$-Benford. (This is what I've called Whittaker's random variable.) For any $c > 1$, define $Y_c \equiv \log_c(X)$. It may then be shown that $Y_c$ is u.d. mod 1 (and hence that $X$ is $c$-Benford) iff $c \leq b$. In summary, $B_X = (1, b]$. **Citations**: [3], [6], [7].

## 4. Digression: Fourier Transforms.

Before going much further, we need to list some facts about Fourier transforms. Let $g()$ denote the pdf of a real valued random variable $Y$. The Fourier transform of $g()$ is the function $\widehat{g}(): \mathbb{R} \to \mathbb{C}$ defined as

$$\widehat{g}(\xi) \equiv \int_{-\infty}^{\infty} e^{-2\pi i \xi y} g(y)\, dy = \mathbb{E}\big(e^{-2\pi i \xi Y}\big) = u(\xi)\ - iv(\xi) \tag{4.1}$$

for all $\xi \in \mathbb{R}$, where

$$\begin{aligned} u(\xi) &\equiv \int_{-\infty}^{\infty} \cos(2\pi\xi y) g(y)\, dy = \mathbb{E}(\cos(2\pi\xi Y)) \quad \text{and} \\ v(\xi) &\equiv \int_{-\infty}^{\infty} \sin(2\pi\xi y) g(y)\, dy = \mathbb{E}(\sin(2\pi\xi Y)). \end{aligned} \tag{4.2}$$

Note that $u()$ is an even function and $v()$ is an odd function, and hence that $\widehat{g}(-\xi) = \overline{\widehat{g}(\xi)}$ where the overbar denotes complex conjugation. Though the Fourier transform $\widehat{g}(\xi)$ is generally complex valued, it is real valued if $g()$ is an even function, i.e. if $Y$ is symmetrically distributed around the origin. Hence, if $g()$ is an even function, then $\widehat{g}()$ is an even function. Finally, note that $\widehat{g}(0) = \int_{-\infty}^{\infty} g(y)\, dy = 1$.

The following fact is very useful. **Proposition 4.1** (shift and scale with random variables). Suppose that $W = \sigma Y + \mu$ where $\sigma > 0$. Suppose that $Y \sim g()$ and let $h()$ denote the pdf of $W$. We may obtain $h()$ from $g()$ and $\widehat{h}()$ from $\widehat{g}()$ as follows:

$$h(w) = \frac{1}{\sigma} g\Big(\frac{w - \mu}{\sigma}\Big) \tag{4.3}$$

(proof left to reader), and

$$\begin{aligned} \widehat{h}(\xi) &= \mathbb{E}\big(e^{-2\pi i \xi W}\big) = \mathbb{E}\big[e^{-2\pi i \xi(\sigma Y + \mu)}\big] \\ &= e^{-2\pi i \xi \mu} \mathbb{E}\big[e^{-2\pi i (\sigma\xi) Y}\big] = e^{-2\pi i \xi \mu} \widehat{g}(\sigma\xi). \end{aligned} \tag{4.4}$$

If $\mu = 0$, eq. (4.4) becomes $\widehat{h}(\xi) = \widehat{g}(\sigma\xi)$.

The Appendix of this paper contains a table of selected Fourier transforms.

## 5. A Framework for Benford Analysis.

Suppose that $X$ is a positive random variable and that $b > 1$. We may wish to know if $X$ is $b$-Benford, and if it's not by how far does it differ from "Benfordness." I call an attempt to answer these and related questions a "Benford analysis." In this section I establish some notation I'll use for a Benford analysis, and give some fundamental results that allow us to proceed.

First, define

$$Y \equiv \log_b(X) = \Lambda_b \ln(X) \quad \text{where} \quad \Lambda_b \equiv \frac{1}{\ln(b)}. \tag{5.1}$$

Next, let

$g()$ denote the pdf of $Y$,
$\widetilde{g}()$ denote the pdf of $\langle Y \rangle$.

Given $\widetilde{g}()$ we may answer the two questions given above. (1) $X$ is $b$-Benford iff $\widetilde{g}(u) = 1$ for almost all $u \in [0,1)$. (2) If $X$ is not $b$-Benford, we may measure its deviation from Benfordness by any measure of the deviation of $\widetilde{g}()$ from a uniform distribution. For example, if $\widetilde{g}()$ is continuous, or if its only discontinuities are "jumps," we could use the infinity norm:

$$\|\widetilde{g} - 1\|_\infty \equiv \sup(|\widetilde{g}(u) - 1| \colon 0 \le u < 1). \tag{5.2}$$

We need a way to find $\widetilde{g}()$ from $g()$. Under a reasonable assumption, it may be shown that

$$\widetilde{g}(u) = \sum_{k \in \mathbb{Z}} g(k + u) \tag{5.3}$$

for all $u \in [0,1)$. The "reasonable assumption" is described in [2]. In this paper we'll just accept eq. (5.3) as given.

Although eq. (5.3) is fundamental for a Benford analysis of $X$, it is not very useful for finding the answers to some analytical questions one may ask. Fourier analysis provides the tools needed to continue the analysis. It may be shown [3] that the Fourier series representation of $\widetilde{g}(u)$ is

$$\widetilde{g}(u) = \sum_{n \in \mathbb{Z}} \widehat{g}(n) e^{2\pi i n u} \quad \text{for all } u \in [0,1). \tag{5.4}$$

At first sight this expression must not seem very useful; the series of real valued functions in eq. (5.3) has been replaced by a series of complex valued functions multiplied by complex coefficients. But $\widehat{g}(0) = 1$, and eq. (5.4) may be written as

$$\widetilde{g}(u) = 1 + \sum_{n \in \mathbb{N}} \left[ \widehat{g}(-n) e^{-2\pi i n u} + \widehat{g}(n) e^{2\pi i n u} \right]. \tag{5.5}$$

As $\widehat{g}(-n)e^{-2\pi inu}$ is the complex conjugate of $\widehat{g}(n)e^{2\pi inu}$, it follows that each term in brackets in eq. (5.5) is real valued. In fact,

$$\widehat{g}(-n)e^{-2\pi inu} + \widehat{g}(n)e^{2\pi inu} = a_n\cos(2\pi nu) + b_n\sin(2\pi nu) \tag{5.6}$$

where

$$\begin{aligned} a_n &= \widehat{g}(-n) + \widehat{g}(n) = 2\int_{-\infty}^{\infty} \cos(2\pi ny)g(y)\,dy, \\ b_n &= -i\widehat{g}(-n) + i\widehat{g}(n) = 2\int_{-\infty}^{\infty} \sin(2\pi ny)g(y)\,dy. \end{aligned} \tag{5.7}$$

Combining eqs. (5.5) and (5.6) yields

$$\widetilde{g}(u) = 1 + \sum_{n\in\mathbb{N}}[a_n\cos(2\pi nu) + b_n\sin(2\pi nu)]. \tag{5.8}$$

In practice, it is often convenient to go one step further and rewrite eq. (5.8) as

$$\widetilde{g}(u) = 1 + \sum_{n\in\mathbb{N}} A_n\cos[2\pi n(u-\theta_n)] \tag{5.9}$$

where $A_n$ satisfies

$$A_n^2 = a_n^2 + b_n^2 \tag{5.10}$$

and $\theta_n$ is any solution to

$$\cos(2\pi n\theta_n) = \frac{a_n}{A_n} \qquad \text{and} \qquad \sin(2\pi n\theta_n) = \frac{b_n}{A_n}. \tag{5.11}$$

The parameters $A_n$ and $\theta_n$ are not uniquely determined by eqs. (5.10) and (5.11), but in practice natural candidates for $A_n$ and $\theta_n$ often present themselves.

I'll call $A_n$ an "amplitude" (though this term generally refers to $|A_n|$) and $\theta_n$ a "phase."

Solving eq. (5.7) for $\widehat{g}(-n)$ and $\widehat{g}(n)$, we find

$$\widehat{g}(n) = \frac{1}{2}(a_n - ib_n), \qquad \widehat{g}(-n) = \frac{1}{2}(a_n + ib_n). \tag{5.12}$$

It follows that

$$A_n^2 = a_n^2 + b_n^2 = 4\widehat{g}(n)\widehat{g}(-n) = 4|\widehat{g}(n)|^2 \quad \Rightarrow \quad |A_n| = 2|\widehat{g}(n)| \tag{5.13}$$

for all $n \in \mathbb{N}$. Hence we've proven the following important fact.

**Proposition 5.1**. The pdf $\widetilde{g}()$ is that of a $U[0,1)$ random variable if and only if $\widehat{g}(n) = 0$ for all $n \in \mathbb{N}$. Equivalently, the pdf $\widetilde{g}()$ is that of a $U[0,1)$ random variable if and only if $A_n = 0$ for all $n \in \mathbb{N}$.

## 6. Base Dependence: Theory.

Suppose we're given a u.d. mod 1 random variable $Y$ with pdf $g()$ and $b > 1$. Then $X \equiv b^Y$ is $b$-Benford. Now let $c > 1$ be another possible base and define $Y_c \equiv \log_c(X)$. Let $g_c()$ and $\widetilde{g}_c()$ denote the pdfs of $Y_c$ and $\langle Y_c \rangle$, respectively, and let $\widehat{g}_c()$ denote the Fourier transform of $g_c()$. My aim in this section is to present tools that allow one to study how $\widetilde{g}_c()$ varies as a function of $c$.

The first thing to observe is that $Y_c$ is proportional to $Y$:

$$Y_c = \frac{\ln(X)}{\ln(c)} = \frac{\ln(b)}{\ln(c)} \cdot \frac{\ln(X)}{\ln(b)} = \rho Y \quad \text{where} \quad \rho \equiv \frac{\ln(b)}{\ln(c)}. \tag{6.1}$$

It then follows from Proposition 4.1 that

$$\widehat{g}_c(\zeta) = \widehat{g}(\rho\zeta) \tag{6.2}$$

for any $\zeta \in \mathbb{R}$.

To use eq. (6.2) we first need to say something about $g()$. I introduced "seed functions" in [2] and showed that every pdf $g()$ of a u.d. mod 1 random variable may be written

$$g(y) = H(y) - H(y-1) \tag{6.3}$$

for every $y \in \mathbb{R}$, where $H()$ is a seed function. Hence

$$\widehat{g}(\xi) = \int_{-\infty}^{\infty} e^{-2\pi i \xi y}[H(y) - H(y-1)]\, dy. \tag{6.4}$$

Under various assumptions about $H()$, we may combine eqs. (6.2) and (6.4) to compute $\widehat{g}_c(n)$ for all $n \in \mathbb{Z}$, and given $\widehat{g}_c(n)$ we may compute $A_n$ and $\theta_n$ in the expression

$$\widehat{g}_c(-n)e^{-2\pi i n u} + \widehat{g}_c(n)e^{2\pi i n u} = A_n \cos[2\pi n(u - \theta_n)]$$

for all $n \in \mathbb{N}$, and thereby derive $\widetilde{g}()$. In this section I'll partially carry out this program for two broad classes of seed functions: (1) Class SF: $H()$ is a **step function**, and (2) Class C: $H()$ is increasing and absolutely **continuous**.

**Example SF 1**. Suppose that $H()$ is the following step function:

$$H(y) = \begin{cases} 0 & \text{if} \quad y < -\frac{1}{2}, \\ 1 & \text{if} \quad y \geq -\frac{1}{2}. \end{cases}$$

This seed function implies that

$$g(y) = \begin{cases} 1 & \text{if} \quad y \in [-\frac{1}{2}, \frac{1}{2}), \\ 0 & \text{otherwise.} \end{cases} \tag{6.5}$$

Hence, from the table of Fourier transforms in the Appendix,

$$\widehat{g}(\xi) = \frac{\sin(\pi\xi)}{\pi\xi} \tag{6.6}$$

(where it is understood that $\widehat{g}(0) = 1$). Combining this with eq. (6.2) yields

$$\widehat{g}_c(n) = \frac{\sin(\pi\rho n)}{\pi\rho n} \tag{6.7}$$

for any $n \neq 0$. From Proposition 5.1 we know that $Y_c$ will be $c$-Benford iff $\widehat{g}_c(n) = 0$ for every integer $n \neq 0$, and from eq. (6.7) it's clear that this happens iff $\rho$ is an integer. But

$$\rho = \frac{\ln(b)}{\ln(c)} = m \quad \Leftrightarrow \quad c = b^{1/m} \tag{6.8}$$

for every $m \in \mathbb{N}$. Hence, $Y_c$ is $c$-Benford iff $c$ is an integral root of $b$. This result agrees with Proposition 3.2.

Certain features of this result are repeated with *every* seed function $H()$ we consider. In particular, we always find that $\widehat{g}_c(n) = 0$ for all $n \neq 0$ whenever $c$ is an integral root of $b$. Also, note that $\widehat{g}_c(n)$ depends on $c$ entirely through the parameter $\rho$.

Equation (6.7) implies that

$$A_n = \frac{2\sin(\pi\rho n)}{\pi\rho n}, \qquad \theta_n = 0 \tag{6.9}$$

for Example SF 1.

**Example SF 2**. To generalize SF 1 slightly, suppose that $H()$ jumps from 0 to 1 at $\mu - \frac{1}{2}$ for some $\mu \in \mathbb{R}$. The pdf $g()$ implied by this seed function is just that given by eq. (6.5) shifted right by $\mu$. From Proposition 4.1 and eq. (6.6) we obtain

$$\widehat{g}(\xi) = e^{-2\pi i\xi\mu}\frac{\sin(\pi\xi)}{\pi\xi}$$

and hence

$$\widehat{g}_c(n) = e^{-2\pi i\rho n\mu}\frac{\sin(\pi\rho n)}{\pi\rho n} \tag{6.10}$$

for any $n \neq 0$. Note that $\widehat{g}_c(n) = 0$ iff $c = b^{1/m}$ for some $m \in \mathbb{N}$. Equation (6.10) implies that

$$A_n = \frac{2\sin(\pi\rho n)}{\pi\rho n}, \qquad \theta_n = \rho\mu \tag{6.11}$$

for this example. The only effect of including $\mu$ is to change the phase. Note that the phase does not depend on $n$.

**Seed functions of Class C**. We now assume that $H()$ is increasing and absolutely continuous. This assumption makes $H()$ mathematically equivalent to the cdf of an absolutely continuous random variable. Under this assumption $H()$ is differentiable almost everywhere and $h(y) \equiv H'(y) \geq 0$. We want to evaluate the integral

$$\widehat{g}(\xi) = \int_{-\infty}^{\infty} e^{-2\pi i\xi y}[H(y) - H(y-1)]\, dy = \mathbb{E}\big(e^{-2\pi i\xi Y}\big).$$

It's clear from the rightmost expression in this equation that $\widehat{g}(0) = 1$. When $\xi \neq 0$, an initial integration by parts yields

$$\widehat{g}(\xi) = \frac{1}{2\pi i\xi}\int_{-\infty}^{\infty} e^{-2\pi i\xi y}[h(y) - h(y-1)]\, dy.$$

Evaluating this integral,

$$\begin{aligned}\widehat{g}(\xi) &= \frac{1}{2\pi i\xi}\big(1 - e^{-2\pi i\xi}\big)\widehat{h}(\xi) \\ &= \frac{e^{-i\pi\xi}}{2\pi i\xi}\big(e^{i\pi\xi} - e^{-i\pi\xi}\big)\widehat{h}(\xi) = \frac{e^{-i\pi\xi}}{\pi\xi}\sin(\pi\xi)\widehat{h}(\xi).\end{aligned}$$

Hence,

$$\widehat{g}_c(n) = \frac{e^{-i\pi\rho n}}{\pi\rho n}\sin(\pi\rho n)\widehat{h}(\rho n) \tag{6.12}$$

for any $n \neq 0$. We see once again that $\widehat{g}_c(n) = 0$ for all $n \in \mathbb{N}$ whenever $c$ is an integral root of $b$. In addition, there's another possibility; $\widehat{g}_c(n) = 0$ for all $n \in \mathbb{N}$ if $\widehat{h}(\rho n) = 0$ for all $n \in \mathbb{N}$. This is essentially the possibility that was exploited in the construction of Whittaker's random variable. We'll return to this point in a moment.

**Example C 1**. Still working with the assumption that $H()$ is increasing and absolutely continuous, we now make the additional assumption that $h()$ is an even function, which implies that $\widehat{h}()$ is an even function. Under these assumptions, eq. (6.12) implies that

$$A_n = \frac{2\sin(\pi\rho n)\widehat{h}(\rho n)}{\pi\rho n}, \qquad \theta_n = \frac{1}{2}\rho. \tag{6.13}$$

**Example C 2**. In Example C 1 we assume that $h()$ is even, so that it's symmetrical around the point $y = 0$. Now assume that $h()$ is symmetrical around the point $y = \mu$ for some $\mu \in \mathbb{R}$. Define $h_0(y) \equiv h(y + \mu)$ so $h_0()$ is an even function. It is easy to show that $\widehat{h}(\xi) = e^{-2\pi i\xi\mu}\widehat{h}_0(\xi)$. Combining this fact with eq. (6.12) yields

$$\begin{aligned}\widehat{g}_c(n) &= \frac{e^{-i\pi\rho n}}{\pi\rho n}\sin(\pi\rho n)e^{-2\pi i\rho n\mu}\widehat{h}_0(\rho n) \\ &= \frac{e^{-2\pi i\rho n(\frac{1}{2}+\mu)}}{\pi\rho n}\sin(\pi\rho n)\widehat{h}_0(\rho n).\end{aligned} \tag{6.14}$$

Equation (6.14) implies that

$$A_n = \frac{2\sin(\pi\rho n)\widehat{h}_0(\rho n)}{\pi\rho n}, \qquad \theta_n = \rho\left(\frac{1}{2} + \mu\right) \tag{6.15}$$

for Example C 2. We observe that the phase depends on $\mu$ and $\rho$, but not on $n$.

Note that $A_n$ and $\theta_n$ depend on $c$ only through $\rho$ in all of these examples. It's useful to keep in mind that $\rho$ depends on $c$ as shown in the following figure (where I've let $b \equiv 16$).

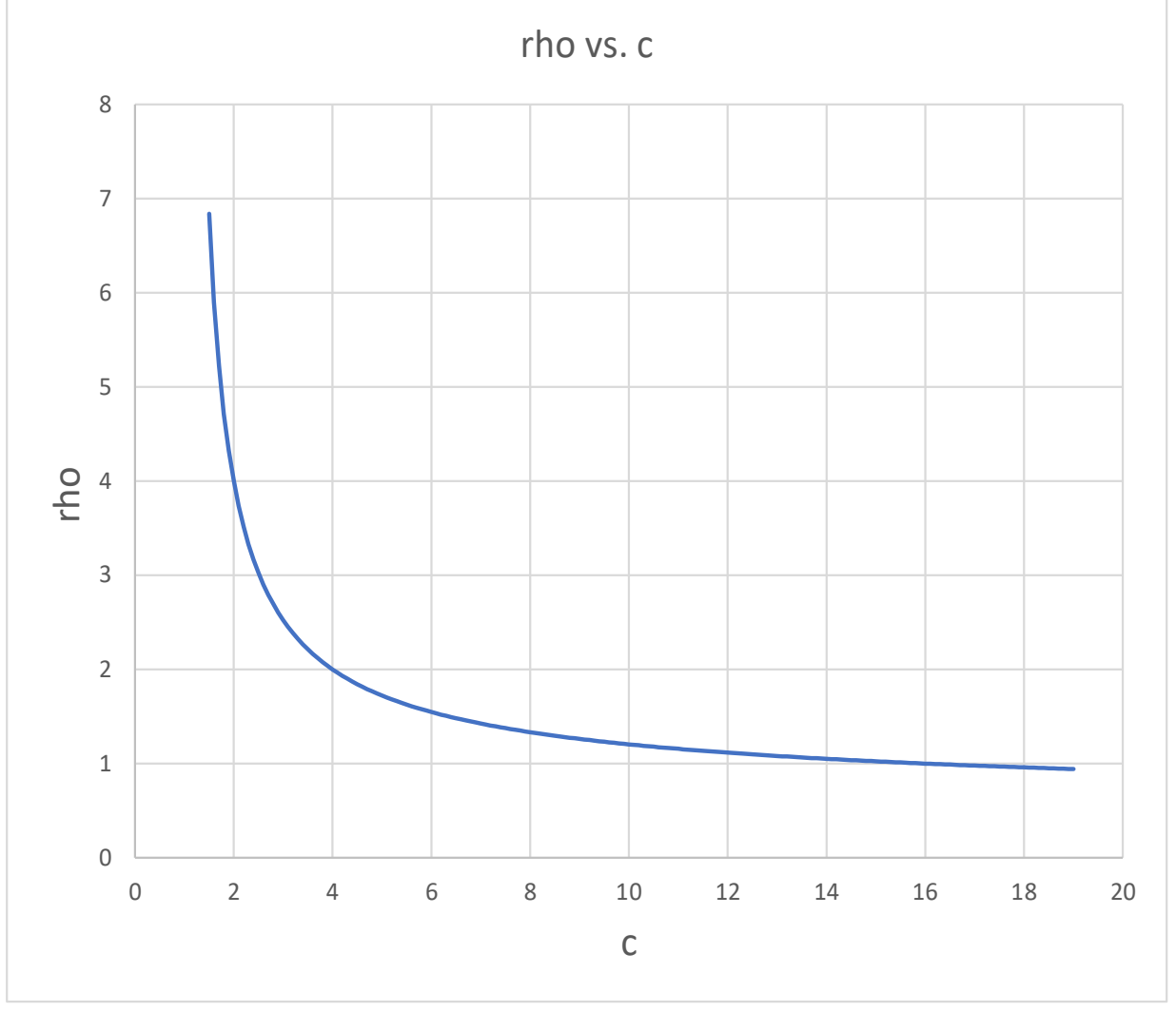

In words, $\rho$ increases from 1 to $\infty$ as $c$ decreases from $b$ towards 1.

## 7. Base Dependence: Examples.

Equations (6.12) and (6.15) provide the scaffolding for the construction of $\widetilde{g}()$, but require insertion of an actual formula for $\widehat{h}()$ in eq. (6.12) or $\widehat{h}_0()$ in eq. (6.15) for completion. This section completes this construction using the table of Fourier transforms in the Appendix.

Every distribution function is a legitimate seed function. Hence every Fourier transform given in the Appendix is a legitimate candidate for $\widehat{h}()$. Moreover, four of the distributions in the Appendix (the normal, Laplace, Cauchy, and logistic) are even functions, and their Fourier transforms are therefore legitimate candidates for $\widehat{h}_0()$. All four of these distributions have fixed variances, however, and it is desirable to append a "scale" parameter $\sigma$ that allows these variances to be adjusted. Proposition 4.1 justifies the following expanded table of Fourier transforms.

| Name | $h_0(y)$ | $\widehat{h}_0(\xi)$ |
|---|---|---|
| $N(0,\sigma)$ | $(2\pi\sigma^2)^{-1/2}e^{-y^2/(2\sigma^2)}$ | $\exp(-2\pi^2\sigma^2\xi^2)$ |
| Laplace$(0,\sigma)$ | $\frac{1}{2\sigma}e^{-\lvert y\rvert/\sigma}$ | $\frac{1}{1+4\pi^2\sigma^2\xi^2}$ |
| Cauchy$(0,\sigma)$ | $\frac{1}{\pi\sigma}\left(1+\frac{y^2}{\sigma^2}\right)^{-1}$ | $e^{-2\pi\sigma\lvert\xi\rvert}$ |
| Logistic$(0,\sigma)$ | $\frac{1}{\sigma}\left(e^{y/(2\sigma)}+e^{-y/(2\sigma)}\right)^{-2}$ | $\frac{2\pi^2\sigma\xi}{\sinh(2\pi^2\sigma\xi)}$ |

**Note**: Among these four distributions, $\sigma$ is the standard deviation of the rescaled random variable *only* for the normal distribution $N(0, \sigma)$.

**Gauss-Benford random variables**. Suppose that $H()$ is the cdf of a $N(\mu, \sigma)$ random variable, i.e., a $N(0, \sigma)$ random variable shifted $\mu$ to the right. I'll call the random variable $X$ implied by this seed function a "Gauss-Benford" random variable. Combining eq. (6.15) with the appropriate entry from this table, we obtain

$$A_n = \frac{2\sin(\pi\rho n)}{\pi\rho n}\exp\left(-2\pi^2\sigma^2\rho^2 n^2\right), \qquad \theta_n = \rho\left(\frac{1}{2} + \mu\right). \tag{7.1}$$

As $\exp(-2\pi^2\sigma^2\rho^2 n^2) > 0$, it follows that $B_X = \left\{b^{1/m} \colon m \in \mathbb{N}\right\}$. Let

$$A_n^* = \frac{2}{\pi\rho n}\exp\left(-2\pi^2\sigma^2\rho^2 n^2\right), \tag{7.2}$$

so $A_n = \sin(\pi\rho n)A_n^*$. Viewed as a function of $n$ or $\rho$, $A_n$ oscillates within an envelope $[-A_n^*, A_n^*]$, and $|A_n| \leq A_n^*$ for all $n,\ \sigma,$ and $\rho$. Asymptotically, letting any of the parameters $n$, $\rho$, or $\sigma \to \infty$ implies that $A_n^* \downarrow 0$. Equation (7.2) implies that $A_1^* > A_2^* > \cdots$. The descent of $A_n^*$ towards zero with increases in $n$, $\rho$, or $\sigma$ is extremely rapid, and $A_1^*$ can be small with even low values of $\rho$ and $\sigma$. For example, letting $\rho = \sigma = 1$ implies that $A_1^* \approx 1.7 \times 10^{-9}$. In this case, the graph of $\widetilde{g}()$ is visually indistinguishable from that of a uniform distribution on $[0, 1)$ and we would have to conclude that $X$ is "effectively" $c$-Benford for all $c \leq b$.

**Laplace-Benford random variables**. Now suppose that $H()$ is the cdf of a Laplace$(\mu, \sigma)$ random variable. I'll call the random variable $X$ implied by this seed function a "Laplace-Benford" random variable. Combining eq. (6.15) with the appropriate entry from the table, we obtain

$$A_n = \frac{2\sin(\pi\rho n)}{\pi\rho n} \cdot \frac{1}{1 + 4\pi^2\sigma^2\rho^2 n^2}, \qquad \theta_n = \rho\left(\frac{1}{2} + \mu\right). \tag{7.3}$$

As $\left(1 + 4\pi^2\sigma^2\rho^2 n^2\right)^{-1} > 0$, it follows that $B_X = \left\{b^{1/m} \colon m \in \mathbb{N}\right\}$. Let

$$A_n^* = \frac{2}{\pi\rho n} \cdot \frac{1}{1 + 4\pi^2\sigma^2\rho^2 n^2}, \tag{7.4}$$

so $A_n = \sin(\pi\rho n)A_n$ and $|A_n| \leq A_n^*$ for all $n,\ \sigma,$ and $\rho$. Asymptotically, letting any of the parameters $n$, $\rho$, or $\sigma \to \infty$ implies that $A_n^* \downarrow 0$. Though the asymptotic limits of eqs. (7.2) and (7.4) are identical, the approach of $A_n^*$ to zero (as $n$, $\rho$, or $\sigma$ increases) is very much slower for eq. (7.4) than it is for eq. (7.2).

**Cauchy-Benford random variables**. Now suppose that $H()$ is the cdf of a Cauchy$(\mu, \sigma)$ random variable. I'll call the random variable $X$ implied by this seed function a "Cauchy-Benford" random variable. Combining eq. (6.15) with the appropriate entry from the table, we obtain

$$A_n = \frac{2\sin(\pi\rho n)}{\pi\rho n}e^{-2\pi\sigma\rho n}, \qquad \theta_n = \rho\left(\frac{1}{2}+\mu\right). \tag{7.5}$$

As $e^{-2\pi\sigma\rho n} > 0$, it follows that $B_X = \{b^{1/m}\colon m \in \mathbb{N}\}$. Let

$$A_n^* = \frac{2}{\pi\rho n}e^{-2\pi\sigma\rho n} \tag{7.6}$$

so $A_n = \sin(\pi\rho n)A_n$. The asymptotic behavior for this $A_n^*$ is identical to that of eq. (7.2) or (7.4). The rate of descent of $A_n^*$ towards zero is intermediate between that of a Gauss-Benford random variable and that of a Laplace-Benford random variable.

**Logistic-Benford random variables**. For our final example of a symmetric seed function, let $H()$ be the cdf of a logistic$(\mu,\sigma)$ random variable. I'll call the random variable $X$ implied by this seed function a "Logistic-Benford" random variable. Combining eq. (6.15) with the appropriate entry from the table, we obtain

$$A_n = \frac{2\sin(\pi\rho n)}{\pi\rho n}\cdot\frac{2\pi^2\sigma\rho n}{\sinh(2\pi^2\sigma\rho n)} = \sin(\pi\rho n)A_n^*, \qquad \theta_n = \rho\left(\frac{1}{2}+\mu\right) \tag{7.7}$$

where

$$A_n^* = \frac{4\pi\sigma}{\sinh(2\pi^2\sigma\rho n)} > 0. \tag{7.8}$$

The asymptotic behavior for this $A_n^*$ is identical to that of the previous three random variables. The rate of convergence of $A_n^*$ to zero is comparable to that of a Cauchy-Benford random variable.

**Gamma-Benford random variables**. Suppose that the seed function $H()$ is the cdf of a Gamma$(\alpha,\beta)$ random variable. I'll call the random variable $X$ implied by this seed function a "Gamma-Benford" random variable. This seed function is increasing and absolutely continuous, but $h()$ is not symmetrically distributed around any point $\mu$, so eq. (6.15) does not apply. Combining eq. (6.12) with the appropriate entry from the table of Fourier transforms found in the Appendix, we obtain

$$\widehat{g}_c(n) = \frac{e^{-i\pi\rho n}}{\pi\rho n}\sin(\pi\rho n)(1+2\pi i\beta\rho n)^{-\alpha} \tag{7.9}$$

for every integer $n \neq 0$. To make headway, define

$$z_n \equiv 1+2\pi i\beta\rho n = 1+iy_n \qquad \text{where} \qquad y_n \equiv 2\pi\beta\rho n, \tag{7.10}$$

and rewrite $z_n$ in polar form, so

$$z_n = r_n e^{i\phi_n} \qquad \text{where} \qquad r_n \equiv \sqrt{1+y_n^2}, \quad \tan(\phi_n) = y_n. \tag{7.11}$$

Hence,

$$\widehat{g}_c(n) = \frac{e^{-i\pi\rho n}}{\pi\rho n}\sin(\pi\rho n)r_n^{-\alpha}e^{-i\alpha\phi_n} = \frac{\sin(\pi\rho n)}{\pi\rho n r_n^{\alpha}}e^{-2\pi i n\theta_n} \tag{7.12}$$

where

$$\theta_n \equiv \frac{1}{2}\rho + \frac{\alpha\phi_n}{2\pi n}. \tag{7.13}$$

Hence,

$$\widehat{g}_c(-n)e^{-2\pi i n u} + \widehat{g}_c(n)e^{2\pi i n u} = \frac{\sin(\pi\rho n)}{\pi\rho n r_n^{\alpha}}\left(e^{-2\pi i n(u-\theta_n)} + e^{2\pi i n(u-\theta_n)}\right)$$

$$= \frac{2\sin(\pi\rho n)}{\pi\rho n r_n^{\alpha}}\cos[2\pi n(u-\theta_n)] = A_n\cos[2\pi n(u-\theta_n)] \tag{7.14}$$

where

$$A_n \equiv \frac{2\sin(\pi\rho n)}{\pi\rho n r_n^{\alpha}}. \tag{7.15}$$

To "compare and contrast" these results with those with symmetric distributions, we make the following observations. (1) The presence of $\sin(\pi\rho n)$ in the numerator of eq. (7.15), combined with $r_n^{-\alpha} > 0$, implies that $A_n = 0$ for all $n \in \mathbb{N}$ if and only if $\rho$ is an integer, i.e., if and only if $c$ is an integral root of $b$. (2) Unlike our earlier results, where the phase $\theta_n$ is given by eq. (6.15) and does not depend on $n$, for a Gamma-Benford random variable the phase is given by eq. (7.13). It's easy to show that $\phi_n \to \frac{1}{2}\pi$ as $n \to \infty$, and hence that $\theta_n \downarrow \frac{1}{2}\rho$. (3) It's easy to show that

$$A_n^* \equiv \frac{2}{\pi\rho n r_n^{\alpha}} \downarrow 0 \qquad \text{as} \qquad \rho n \to \infty.$$

**Whittaker-Benford random variables**. For our final example, we return to eq. (6.12)

$$\widehat{g}_c(n) = \frac{e^{-i\pi\rho n}}{\pi\rho n}\sin(\pi\rho n)\widehat{h}(\rho n), \tag{6.12}$$

which holds for all seed functions in Class C. All of our previous examples have made use of the fact that $\sin(\pi\rho n) = 0$ for all $n \in \mathbb{N}$ whenever $\rho$ is an integer. We now consider another possibility: $\widehat{g}_c(n) = 0$ for all $n \in \mathbb{N}$ if $\widehat{h}(\rho n) = 0$ for all $n \in \mathbb{N}$. I'll say that a $b$-Benford random variable $X$ satisfying this condition is a "Whittaker-Benford" random variable. The key here is to find $\widehat{h}()$ with bounded support, and the simplest such $\widehat{h}()$ is triangular:

$$\widehat{h}(\xi) = \max\left(0, 1 - \frac{|\xi|}{\gamma}\right), \tag{7.16}$$

where $\gamma > 0$. With this $\widehat{h}()$ it's clear that $\widehat{h}(\rho n) = 0$ for all $n \in \mathbb{N}$ if $\rho \geq \gamma$. Note that $\rho \geq \gamma \Leftrightarrow c \leq b^{1/\gamma}$. Therefore, the Benford spectrum $B_X$ of a Whittaker-Benford random variable $X$ with $\widehat{h}()$ given by eq. (7.16) has two (overlapping) components: $B_X = B_X^d \cup B_X^c$ where

$$\begin{aligned} B_X^d &\equiv \left\{ b^{1/m} \colon m \in \mathbb{N} \right\}, \\ B_X^c &\equiv (1, b^{1/\gamma}]. \end{aligned} \tag{7.17}$$

(The superscript $d$ stands for "discrete," and the superscript $c$ stands for "continuous.") If $\gamma \leq 1$, then $B_X^d \subset B_X^c$. For example, if $\gamma = \frac{1}{2}$ then $B_X = B_X^c = (1, b^2]$. On the other hand, if $\gamma > 1$, then $B_X^c = (1, b^{1/\gamma}] \subset (1, b]$, so $B_X$ equals the disjoint union of the discrete set $(B_X^d - B_X^c)$ and the continuous set $(B_X^c)$.

The function $h()$ that yields $\widehat{h}()$ given by eq. (7.16) is

$$h(y) = \frac{1 - \cos(2\pi\gamma y)}{2\gamma\pi^2 y^2}. \tag{7.18}$$

## 8. On "Base-Invariant Significant Digits".

I wish to acknowledge that I first encountered many of the ideas discussed in this section in Michal Wójcik's admirable paper [7]. To the best of my knowledge, however, Proposition 8.4 is entirely my own.

**Proposition 8.1**. If $X$ is $b$-Benford, then $X^n$ is $b$-Benford for any $n \in \mathbb{N}$. **Proof**. As $X$ is $b$-Benford, $X = b^Y$ where $Y$ is u.d. mod 1. Hence $X^n = b^{nY}$. But $nY$ is u.d. mod 1 by Proposition 3.1. Therefore $X^n$ is $b$-Benford.

**Corollary 1**. As $X^n = (b^n)^Y$, it follows that $X^n$ is $b^n$-Benford if $X$ is $b$-Benford.

**Corollary 2**. If $X$ is $b$-Benford, then $S_b(X) \sim S_b(X^n)$ for any $n \in \mathbb{N}$. This follows from Definition 2.1.

One may wonder if the converse of Corollary 2, namely

$$\text{if } S_b(X) \sim S_b(X^n) \text{ for all } n \in \mathbb{N}, \text{ then } X \text{ is } b\text{-Benford},$$

is true. The answer is "no." Here's a counterexample. If $X \equiv 1$, then $S_b(X) \sim S_b(X^n)$ for all $n \in \mathbb{N}$, but $X$ is not $b$-Benford. In fact, any $X$ of the form $b^m$ where $m \in \mathbb{Z}$ is a counterexample, as $S_b(X) = 1 = S_b(X^n)$. However, we may show the following:

**Proposition 8.2.** If $S_b(X) \sim S_b(X^n)$ for all $n \in \mathbb{N}$, then either $X$ is $b$-Benford, or $S_b(X) = 1$. We'll provide a proof in a moment.

Let

$$S_b(X) \sim S_b(X^n) \text{ for all } n \in \mathbb{N} \tag{8.1}$$

be called Wójcik's condition.

Here's another way to state Proposition 8.2. Call this **Proposition 8.2a**. (It is Wójcik's Theorem 19.) $X$ satisfies Wójcik's condition if and only if the cdf of $S_b(X)$ is given by

$$\Pr(S_b(X) \le s) = q + (1-q)\log_b(s) \tag{8.2}$$

for some $q \in [0,1]$ and for all $s \in [1, b)$.

To prove Proposition 8.2 or 8.2a, we first massage Wójcik's condition into an alternative form. Let $X$ be a positive random variable and define $Y \equiv \log_b(X)$. For all $n \in \mathbb{N}$,

$$\begin{aligned} S_b(X) \sim S_b(X^n) \quad &\Leftrightarrow \quad \log_b(S_b(X)) \sim \log_b(S_b(X^n)) \\ &\Leftrightarrow \quad \langle \log_b(X) \rangle \sim \langle \log_b(X^n) \rangle = \langle n \log_b(X) \rangle \\ &\Leftrightarrow \quad \langle Y \rangle \sim \langle nY \rangle = \langle n \langle Y \rangle \rangle \end{aligned} \tag{8.3}$$

where the last equality follows from the identity $\langle ny \rangle = \langle n\langle y \rangle + n\lfloor y \rfloor \rangle = \langle n \langle y \rangle \rangle$ for any $y \in \mathbb{R}$.

Berger and Hill ([4], Lemma 5.15, page 77) show the following. **Proposition 8.3**. $\langle Y \rangle \sim \langle n \langle Y \rangle \rangle$ for all $n \in \mathbb{N}$ if and only if

$$\Pr(\langle Y \rangle \le u) = q + (1-q)u \qquad \text{for all } \ u \in [0,1) \tag{8.4}$$

for some $q \in [0,1]$.

Propositions 8.2 and 8.2a are straightforward corollaries of Proposition 8.3.

I bring these facts to the reader's attention because Wójcik's condition is effectively equivalent to Berger and Hill's notion of "base-invariant significant digits" and sheds some light on this notion. (I say "effectively equivalent" as Berger and Hill's concept applies to a probability measure $P$, whereas Wójcik's condition applies to a random variable $X$.)

Here's Berger and Hill's definition (**Definition 5.10**, page 75). Let $\mathcal{A} \supseteq \mathcal{S}$ be a $\sigma$-algebra on $\mathbb{R}^+$. A probability measure $P$ on $(\mathbb{R}^+, \mathcal{A})$ has *base-invariant significant digits* if $P(A) = P(A^{1/n})$ for all $A \in \mathcal{S}$ and all $n \in \mathbb{N}$.

Here's a guide to the symbols used in this definition. (1) $\mathcal{S}$ is the $\sigma$-algebra generated by the significand function $S_b()$. (2) $\mathbb{R}^+ \equiv (0, \infty)$, the set of strictly positive real numbers. (3) For any $A \subseteq \mathbb{R}^+$ and $n \in \mathbb{N}$, $A^{1/n} \equiv \{x > 0 \colon x^n \in A\}$. Also, it's useful at this point to introduce one more bit of non-standard notation used by Berger and Hill: for every $x \in \mathbb{R}$ and every set $C \subseteq \mathbb{R}$, let $xC \equiv \{xc \colon c \in C\}$.

**Proposition 8.4**. Suppose that $X \sim (\mathbb{R}^+, \mathcal{B}(\mathbb{R}^+), P)$. Then $S_b(X) \sim S_b(X^n)$ for all $n \in \mathbb{N}$ if and only if $P$ has base-invariant significant digits.

**Proof**. We begin by proving that Wójcik's condition holds whenever $P$ has base-invariant significant digits. Suppose that $A \in \mathcal{S}$. From the definition of $A^{1/n}$, we have $X \in A^{1/n} \Leftrightarrow X^n \in A$. Hence,

$$P(A^{1/n}) = \Pr(X \in A^{1/n}) = \Pr(X^n \in A). \tag{8.5}$$

If $P$ has base-invariant significant digits, then

$$P(A^{1/n}) = P(A) = \Pr(X \in A). \tag{8.6}$$

Combining eqs. (8.5) and (8.6), we see that

$$\Pr(X \in A) = \Pr(X^n \in A) \tag{8.7}$$

whenever $P$ has base-invariant significant digits. As $A \in \mathcal{S}$ there exists a set $A_0 \in \mathcal{B}[1, b)$ such that

$$A = \bigcup_{k \in \mathbb{Z}} b^k A_0. \tag{8.8}$$

In fact, $A_0 = S_b(A) \equiv \{S_b(x) \colon x \in A\}$. Hence

$$\begin{aligned} X \in A &\Leftrightarrow S_b(X) \in A_0, \\ X^n \in A &\Leftrightarrow S_b(X^n) \in A_0. \end{aligned} \tag{8.9}$$

Combining eqs. (8.7) and (8.9), we conclude that

$$\Pr(S_b(X) \in A_0) = \Pr(S_b(X^n) \in A_0). \tag{8.10}$$

As this equation holds for every $A_0 \in \mathcal{B}[1, b)$, we conclude that $S_b(X) \sim S_b(X^n)$ whenever $X \sim (\mathbb{R}^+, \mathcal{B}(\mathbb{R}^+), P)$ and $P$ has base-invariant significant digits.

To prove that Wójcik's condition implies that $P$ has base-invariant significant digits, we essentially reverse this chain of logic. Wójcik's condition implies eq. (8.10) for any $A_0 \in \mathcal{B}[1, b)$, which in turn implies eq. (8.7) for $A$ given by eq. (8.8). But $\Pr(X \in A) = P(A)$ and $\Pr(X^n \in A) = P\left(A^{1/n}\right)$, so eq. (8.7) implies that $P(A) = P\left(A^{1/n}\right)$. As $A_0$ was an arbitrary element of $\mathcal{B}[1, b)$, the equation $P(A) = P\left(A^{1/n}\right)$ holds for $A$, an arbitrary element of $\mathcal{S}$, and the proof is complete.

Berger and Hill state the following theorem (Theorem 5.13, page 76). A probability measure $P$ on $(\mathbb{R}^+, \mathcal{A})$ with $\mathcal{A} \supseteq \mathcal{S}$ has base-invariant significant digits if and only if, for some $q \in [0, 1]$,

$$P(A) = q\delta_1(A) + (1-q)\mathbb{B}(A) \quad \text{for every } \ A \in \mathcal{S}.$$

(The meaning of the "Dirac measure" $\delta_1()$ is given on page 22 of their book, and the meaning of the "Benford measure" $\mathbb{B}()$ is given on page 32.)

In the light of Proposition 8.4, it can be seen that Berger and Hill's Theorem 5.13 is equivalent to Proposition 8.2a given above.

I conclude this section with a personal opinion about Berger and Hill's exposition. I think that the terminology "base-invariant" they chose for their concept is a misnomer. There is only one base ($b$) in the definition, and their concept of "base-invariant" significant digits tells us nothing about the Benford properties of alternative bases for a $b$-Benford random variable. Hence, the label "base-invariant" they chose for their concept seems misleading and I believe they really should give it a different name.

## 9. Conclusion and Prospect.

Let $Y$ be a u.d. mod 1 random variable with pdf $g()$, let $b > 1$, and define $X \equiv b^Y$, so $X$ is $b$-Benford. Without loss of generality we may assume that

$$g(y) = H(y) - H(y-1) \quad \text{for any } y \in \mathbb{R},$$

where $H()$ is a seed function. Let $c > 1$. In principle, the machinery introduced in Section 6 allows one to investigate the dependence of the distribution of $\langle \log_c(X) \rangle$ on $c$. In practice, I've carried out this investigation only for a limited number of seed functions. Specifically, I've restricted attention to seed functions in the following two classes.

(1) Step functions that jump from 0 to 1 in a single step.

(2) Increasing functions that are absolutely continuous.

What other sort of seed functions are there? Here's a brief list.

(3) Step functions that increase from 0 to 1 at a finite or countably infinite number of "points of jump."

(4) Convex combinations of seed functions in classes (2) and (3).

(5) "Singular" distribution functions. These functions are increasing and continuous, but not absolutely continuous. The Cantor function is the best known example.

(6) Seed functions satisfy a condition I call "unit interval increasing." Every increasing function is unit interval increasing, but not conversely. That is, a function $H()$ may be unit interval increasing, but not everywhere increasing. Several examples of such seed functions are given in [2].

My intuition suggests that seed functions of types (3) and (4) will offer no additional conceptual difficulties, though they will certainly complicate the algebra. I'll leave the investigation of seed functions of classes (5) and (6) to the reader.

With $X$ and $c$ defined as above, let $\widetilde{g}_c()$ denote the pdf of $\langle \log_c(X) \rangle$. If $X$ is $c$-Benford, and if $\widetilde{g}_c()$ is continuous or has only "jump" discontinuities, then

$$\|\widetilde{g}_c - 1\|_\infty = 0. \tag{9.1}$$

Hence, $c$ is in the Benford spectrum $B_X$ iff eq. (9.1) is satisfied. For almost all random variables $X$, the Benford spectrum $B_X$ is empty. We might want to say that $X$ is "effectively" $c$-Benford if

$$\|\widetilde{g}_c - 1\|_\infty < \epsilon \tag{9.2}$$

for some small number $\epsilon$. If we define the "effective" Benford spectrum of $X$ to be the set

$$B_{X,\epsilon} \equiv \{c > 1 \colon \|\widetilde{g}_c - 1\|_\infty < \epsilon\}, \tag{9.3}$$

then $B_X \subseteq B_{X,\epsilon}$. In general, I suggest, the effective spectrum will be a much larger set than the spectrum.

The machinery described in Section 5 to carry out a "Benford analysis" helps us determine whether or not the criterion of eq. (9.2) is satisfied. In Section 7 I suggested that a Gauss-Benford random variable should be regarded as effectively $c$-Benford if the product $\rho\sigma$ is large enough. In [3] I suggested that a lognormal random variable, which is not $b$-Benford for any $b$, should be regarded as effectively $c$-Benford if $\Lambda_c\sigma$ is large enough.

I leave further investigation of effectively Benford random variables to the reader.

**References.**

**Appendix: A Small Table of Fourier Transforms**

Feller [5] gives a table of characteristic functions of selected probability density functions. I've adapted his table to give the Fourier transforms of 8 of his 10 densities, and added a row for an additional pdf (the logistic).

| No. | Name | Density $g(x)$ | Interval | Fourier Transform $\widehat{g}(\xi)$ |
|---|---|---|---|---|
| 1 | $N(0,1)$ | $(2\pi)^{-1/2}e^{-x^2/2}$ | $\mathbb{R}$ | $\exp(-2\pi^2\xi^2)$ |
| 2 | $U[-a,a]$ | $1/2a$ | $[-a,a]$ | $\frac{\sin(2\pi a\xi)}{2\pi a\xi}$ |
| 3 | $U[0,a]$ | $1/a$ | $[0,a]$ | $\frac{1-\exp(-2\pi ia\xi)}{2\pi ia\xi}$ |
| 4 | Triangular | $\frac{1}{a}\left(1-\frac{\lvert x\rvert}{a}\right)$ | $\lvert x\rvert \leq a$ | $\frac{1-\cos(2\pi a\xi)}{2\pi^2a^2\xi^2}$ |
| 5 | Dual of 4 | $\frac{1-\cos(2\pi ax)}{2a\pi^2x^2}$ | $\mathbb{R}$ | $\max\left(0,1-\frac{\lvert\xi\rvert}{a}\right)$ |
| 6 | $\Gamma(\alpha,\beta)$ | $\frac{1}{\Gamma(\alpha)}x^{\alpha-1}e^{-x}$ | $x>0$ | $(1+2\pi i\beta\xi)^{-\alpha}$ |
| 7 | Laplace$(0,1)$ | $\frac{1}{2}e^{-\lvert x\rvert}$ | $\mathbb{R}$ | $\frac{1}{1+4\pi^2\xi^2}$ |
| 8 | Cauchy$(0,1)$ | $\frac{1}{\pi}\frac{1}{1+x^2}$ | $\mathbb{R}$ | $e^{-2\pi\lvert\xi\rvert}$ |
| 9 | Logistic$(0,1)$ | $\left(e^{x/2}+e^{-x/2}\right)^{-2}$ | $\mathbb{R}$ | $\frac{2\pi^2\xi}{\sinh(2\pi^2\xi)}$ |